\documentclass{article}
\usepackage[utf8]{inputenc}
\DeclareMathAlphabet{\mathpzc}{OT1}{pzc}{m}{it}

\usepackage{amsmath} 
\DeclareMathOperator{\sinc}{sinc}

\usepackage{amsfonts}
\usepackage{titlesec}
\usepackage{dsfont}
\usepackage[utf8]{inputenc}
\usepackage[english]{babel}
\usepackage[T1]{fontenc}
\usepackage{physics}
\usepackage{authblk}
\usepackage{tikz}
\usetikzlibrary{calc,decorations.markings}
\usepackage{subcaption}
\usepackage{graphicx}
\usepackage{multicol}
\usepackage{amsthm}
\usepackage{float}
\newtheorem{theorem}{Theorem}[section]

\newtheorem{corollary}{Corollary}[theorem]
\newtheorem{lemma}[theorem]{Lemma}
\newtheorem*{remark}{Remark}
\usepackage{anyfontsize}
\newmuskip\pFqmuskip
\newcommand*\pFq[6][8]{%
  \begingroup 
  \pFqmuskip=#1mu\relax
  \mathcode`\,=\string"8000
  \begingroup\lccode`\~=`\,
  \lowercase{\endgroup\let~}\pFqcomma
  {}_{#2}F_{#3}{\left[\genfrac..{0pt}{}{#4}{#5};#6\right]}%
  \endgroup
}
\newcommand{\pFqcomma}{\mskip\pFqmuskip}
\makeatletter
\newsavebox{\@brx}
\newcommand{\llangle}[1][]{\savebox{\@brx}{\(\m@th{#1\langle}\)}%
  \mathopen{\copy\@brx\kern-0.5\wd\@brx\usebox{\@brx}}}
\newcommand{\rrangle}[1][]{\savebox{\@brx}{\(\m@th{#1\rangle}\)}%
  \mathclose{\copy\@brx\kern-0.5\wd\@brx\usebox{\@brx}}}
\makeatother

\begin{document}

\title{An explanation of the commuting operator “miracle” in time and band limiting}
\author[1]{Pierre-Antoine Bernard}
\author[2]{Nicolas Crampé}
\author[1,3]{Luc Vinet}
\affil[1]{Centre de recherches mathématiques, Université de Montréal, P.O. Box 6128, Centre-ville
Station, Montréal (Québec), H3C 3J7, Canada,}
\affil[2]{Institut Denis-Poisson CNRS/UMR 7013 - Université de Tours - Université d’Orléans, Parc de
Grandmont, 37200 Tours, France.}
\affil[3]{IVADO, 6666 Rue Saint-Urbain, Montréal (Québec), H2S 3H1, Canada}

\maketitle
\begin{abstract}
Time and band limiting operators are expressed as functions of the confluent Heun operator arising in the spheroidal wave equation. Explicit formulas are obtained when the bandwidth parameter is either small or large and results on the complete Fourier transform are recovered. 
\end{abstract}

\section{Introduction}
In a famous series of papers on the time and band limiting of functions \cite{PSWF2,PSWF3,PSWF4,PSWF5,PSWF1}, Slepian, Pollack and Landau made the surprising observation that a second order linear differential operator arising in the confluent Heun equation,
\begin{align}
    T = (1-x^2) \pdv[2]{}{x} - 2x \pdv{}{x} - c^2 x^2,
    \label{heun}
\end{align}
commutes with the finite Fourier transform $\mathcal{F}_c$ \cite{popolit,wang2017review}:
\begin{align}
    \mathcal{F}_c[\phi](x) = \int_{-1}^{1} e^{icxt} \phi(t)dt,
    \label{fifit}
\end{align}
where the bandwidth parameter $c$ is an arbitrary positive number. The eigenfunctions of $T$ and $\mathcal{F}_c$ were further recognized to be \textit{prolate spheroidal wave functions}, which appear in solutions of the Helmholtz equation in appropriate coordinates. Naturally, these functions were also found to diagonalize the integral operator with $\sinc$ kernel $\mathcal{Q}_c = \frac{2\pi}{c}\mathcal{F}_c^* \circ \mathcal{F}_c$,
\begin{align}
    \mathcal{Q}_c[\phi](x) = \int_{-1}^{1}\frac{\sin(c(x-t))}{\pi(x-t)}\phi(t)dt,
    \label{inin}
\end{align}
and were used to derive asymptotic expressions for its spectrum \cite{slepian1965some, des1973asymptotic}. Since then, many fields have benefited from these results. Applications have in particular been made in limited angle tomography \cite{davison1983ill,grunbaum1982limited}, random matrix theory \cite{mehta2004random,des1973asymptotic}, signal processing and in the study of entanglement in fermionic systems \cite{eisler2013free}. The unexpected discovery of this commuting operator raised the following question: what is behind this “miracle” or what is the nature of the relation between the Heun operator $T$ and the finite Fourier transform $\mathcal{F}_c$? 

Recently, an answer explaining the existence a commuting second order differential operator was presented \cite{grunbaum2018algebraic}. By relating $\mathcal{F}_c$ to a certain type of bispectral problem, it was shown that $T$ could be constructed as a special case of an algebraic Heun operator. Furthermore, this framework was applied to other settings where a second order differential operator commutes with an integral one and to cases where a full matrix commutes with a tridiagonal one. 

To understand how the confluent Heun operator $T$ and the finite Fourier transform $\mathcal{F}_c$ are related, an alternative avenue would be to express one as a function of the other. This has been carried out in the case of the complete Fourier transform $\mathcal{F}$,
\begin{align}
    \mathcal{F}[\phi](x) = \frac{1}{\sqrt{2\pi}} \int_{-\infty}^{\infty} e^{ixt} \phi(t) dt 
\end{align}
which commutes with the operator
\begin{align}
    \mathcal{H} = \frac{1}{2}\left(\pdv[2]{}{x} - x^2 - 1\right).
    \label{hamilher}
\end{align}
This last property becomes manifest upon observing that $\mathcal{F}$ can be expressed as the following exponential of $\mathcal{H}$:
\begin{align}
    \mathcal{F} =  e^{-i\frac{\pi}{2}\mathcal{H}}.
    \label{cfc}
\end{align}
Our objective is thus to generalize this formula and to obtain an analogue for $\mathcal{F}_c$ and $T$. The paper is divided in three parts. In section \ref{s2},  we introduce a family of operators which are functions of $T$ and for which the action on the space of square-integrable functions $L^2[-1,1]$ is easy to derive. In section \ref{s3}, we obtain formulas for $\mathcal{F}_c$ and $\mathcal{Q}_c$ in terms of $T$. In section \ref{s4}, we consider the limits $c\rightarrow 0$ and $c \rightarrow \infty$ and recover equation \eqref{cfc}. 

\section{The operators $U(\xi; T)$}
\label{s2}
To express $\mathcal{F}_c$ as a function of $T$, let us start by constructing a set of operators $\{U(\xi; T)\}_{\xi \in \Xi \subset \mathbb{R}}$ which are functions of $T$ and for which the action on functions $f \in L^2[-1,1]$ of the variable $x$ is easy to derive. Consider the following equation:
\begin{equation}
    \begin{split}
        \Big[ (1-x^2) \pdv[2]{}{x} - 2x \pdv{}{x} & - c^2 x^2 \Big] f(x,y) =\\
        &\Big[ (1-y^2) \pdv[2]{}{y} - 2y \pdv{}{y} - c^2 y^2 \Big]f(x,y),
    \end{split}
    \label{main}
\end{equation}
or equivalently
\begin{align}
    T_x f(x,y) = T_y f(x,y),
\end{align}
where $T_x$ and $T_y$ refer to the Heun operator defined in \eqref{heun} acting on the variable $x$ and $y$ respectively. It is interesting to note that equation \eqref{main} (restricted to $x \in [-1,1]$ and $y > 1$) arises from the Helmholtz equation in prolate spheroidal coordinates when a cylindrical symmetry is assumed. 

Rearranging equation \eqref{main}, one finds that the first derivative of $f(x,y)$ with respect to $y$ can be expressed as
\begin{align}
           \pdv{}{y}f(x,y) = \frac{1}{2y}\Big( (1-y^2) \pdv[2]{}{y} - T_x - c^2 y^2 \Big)f(x,y).
           \label{devy}
\end{align}
In particular, the first derivative evaluated at the regular singular points $y = \pm 1$ can be expressed as a linear function of $T_x$ acting on $f(x,\pm 1)$:
\begin{align}
           \pdv{}{y}f(x,y)\Big|_{y = \pm 1} = \mp \frac{ 1}{2}\Big( T_x + c^2 \Big)f(x,\pm 1).
           \label{devy}
\end{align}
Furthermore, the application of $\pdv[k]{}{y}$ on \eqref{main} gives
\begin{equation}
    \begin{split}
        \pdv[k+1]{}{y}f(x,y) &= \bigg[ \frac{(1-y^2)}{2 y (k+1)}\pdv[k+2]{}{y} - \frac{\left( T_x + c^2 y^2 + k(k+1)\right)}{2 y (k+1)} \pdv[k]{}{y} \\
        & \quad \quad  - \frac{c^2 y k }{y (k+1)} \pdv[k-1]{}{y} - \frac{c^2 k(k-1)}{2 y (k+1)}\pdv[k-2]{}{y}\bigg] f(x,y),
        \label{4t}
    \end{split}
\end{equation}
which can be used to obtain the following lemma:
\begin{lemma}
Let $f(x,y)$ be a solution of equation \eqref{main}. The $k^{\text{th}}$ derivative of $f(x,y)$ with respect to $y$ evaluated at $y= \pm1$ can be expressed in terms of a polynomial $U_k$ of $T_x$ acting on $f(x,\pm 1)$, i.e.
\begin{align}
     \pdv[k]{}{y}& f(x,y)\Big|_{y = \pm 1} = (\mp 1)^k U_{k}(T_x) f(x,\pm 1).
     \label{ddv}
\end{align}
The polynomials $U_k$ are given by the following four-term recurrence relation:
\begin{equation}
    \begin{split}
            U_{k+1}(T_x) &= \frac{(T_x + c^2 + k(k+1))}{2(k + 1)}U_k(T_x) \\& \quad \quad  - \frac{c^2 k}{k + 1}U_{k-1}(T_x) + \frac{ c^2 k(k-1)}{2(k + 1)} U_{k-2}(T_x),
    \end{split}
    \label{rel2}
\end{equation}
and the initial condition $U_0(T_x) = 1$.
\end{lemma}
\begin{remark}
Given the four-term recurrence relation \eqref{rel2}, the Theorem 3.2 in \cite{van1987vector} can be applied and one concludes that the polynomials $U_k$ form a family of $2$-orthogonal polynomials.
\end{remark}
\noindent Using the polynomials $U_k$, we can define the operators
\begin{align}
    U(\xi ; T_x) \equiv \sum_{k =0}^\infty  \frac{\xi^k U_k(T_x)}{k!},
    \label{defU}
\end{align}
where $\xi \in ]-2,2[$. They are well defined since the four-term recurrence relation \eqref{4t} yields
\begin{align}
     \frac{\xi^k U_k}{k!} \sim \left(\frac{\xi}{2}\right)^k, \quad \text{as} \quad k \rightarrow \infty,
\end{align}
and thus \eqref{defU} converges as long as $\xi  \in ]-2,2[$.  By equation \eqref{ddv} the operators $U(\xi ; T_x)$ also verify
\begin{equation}
    \begin{split}
          U(\xi ; T_x)[f](x,\pm 1) &= \sum_{k =0}^\infty  \frac{\xi^k }{k!}U_k(T_x)[f](x,\pm 1)\\
          & = \sum_{k =0}^\infty  \frac{(\mp \xi)^k }{k!}\pdv[k]{}{y}f(x,y)\Big|_{y = \pm 1}\\
          & = f(x,\pm(1-\xi)).
    \end{split}
\end{equation}
In other words, they can be interpreted as translation operators on the variable $y$. The next theorem follows.
\begin{theorem}
Let $f_0 \in L^2[-1,1]$, $\xi \in ]-2,2[$ and $T_x$ be the second order linear differential operator defined in \eqref{heun}. Then,
\begin{align}
    U(\xi;T_x)[f_0](x) = f(x, \pm (1 - \xi))
\end{align}
where $U(\xi;T_x)$ is the operator defined in \eqref{defU} and $f$ is the solution of 
\begin{align}
    T_x f(x,y) = T_y f(x,y),
\end{align}
which verifies the boundary condition
\begin{align}
    f(x,\pm 1) = f_0(x).
\end{align} 
\label{theo1}
\end{theorem}
\noindent By construction and from the previous theorem, one also deduces the following:
\begin{corollary}
Let $\lambda \in \mathbb{R}$. It is observed that
\begin{align}
    U(y + 1; \lambda) = \sum_{k =0}^\infty  \frac{(y+1)^k U_k(\lambda)}{k!}, \
\end{align}
 with the polynomials $U_k$ given by the four-term recurrence relation \eqref{rel2}, verifies the confluent Heun equation
\begin{align}
   \left((1-y^2) \pdv[2]{}{y} - 2y \pdv{}{y} - c^2 y^2 - \lambda \right) U(y + 1; \lambda) = 0,
\end{align}
for all $y \in ]-3,1[$, and the boundary condition
\begin{align}
    U(0; \lambda) = 1.
\end{align}
\label{cococo}
\end{corollary}
Next, we shall look for linear combinations of $U(\xi ; T_x)$ which reproduce the action of $\mathcal{F}_c$ and $\mathcal{Q}_c$ on $L^2[-1,1]$.

\section{Formulas for $\mathcal{F}_c$ and $\mathcal{Q}_c$}
\label{s3}
We want to find $\alpha(\xi)$ and $\beta(\xi)$ such that
\begin{align}
    \mathcal{F}_c = \int \alpha(\xi) U(\xi;T_x) d\xi \quad \text{and} \quad \mathcal{Q}_c = \int \beta(\xi) U(\xi;T_x) d\xi.
    \label{conme}
\end{align}
Since $U(\xi; T_x)$ is a function of $T_x$ for all $\xi \in ]-2,2[$, this is sufficient to express $\mathcal{F}_c$ and $\mathcal{Q}_c$ as functions of the Heun operator $T_x$. For \eqref{conme} to be verified, both sides of each equation must have the same diagonal action on the basis of $L^2[-1,1]$ given by the prolate spheroidal wave functions $\psi^c_n(x)$, $n \in \mathbb{N}$, i.e.
\begin{align}
    \mathcal{F}_c[\psi_n^c] = \int \alpha(\xi) U(\xi;T_x)[\psi_n^c] d\xi \quad \text{and} \quad \mathcal{Q}_c[\psi_n^c] = \int \beta(\xi) U(\xi;T_x)[\psi_n^c] d\xi.
\end{align}
Let us recall some properties of these functions.

\subsection{The prolate spheroidal wave functions}
The properties discussed in this subsection can be found in \cite{popolit, osipov2013prolate, wang2017review}. First, we note that the prolate spheroidal wave functions $\psi_n^c$, $n\in \mathbb{N}$, satisfy the following eigenvalue equation
\begin{align}
    T_x \psi^c_n(x) = -\chi_n(c) \psi^c_n(x),
\end{align}
and give a basis of $L^2[-1,1]$. The eigenvalues $\chi_n(c)$ are positive and ordered such that for all $c > 0$
\begin{align}
    \chi_n(c) < \chi_{n+1}(c), \quad \forall n \in \mathbb{N}.
\end{align}
Next, these functions also diagonalize the finite Fourier transform and the $\sinc$ kernel integral operator defined respectively in equation \eqref{fifit} and \eqref{inin}:
\begin{align}
    \mathcal{F}_c[\psi^c_n](x) = i^n \lambda_n(c)\psi^c_n(x), \quad \mathcal{Q}_c[\psi^c_n](x) = \mu_n(c) \psi^c_n(x),
\end{align}
where
\begin{align}
    \mu_n = \frac{c}{2\pi}|\lambda_n(c)|^2. 
\end{align}
Finally, it is also interesting to note that 
\begin{align}
   \mathcal{R}\psi_n^c(x) =  \psi_n^c(-x) = (-1)^n \psi_n^c(x),
\end{align}
where $\mathcal{R}$ refers to the reflection operator acting on functions of the variable $x$. In other words, these functions are even for $n$ even and odd for $n$ odd.

\subsection{The finite Fourier transform}
We look for $\alpha(\xi)$ such that
\begin{align}
   \left( \int \alpha(\xi) U(\xi ; T_x) d\xi\right)[\psi_n^c] (x) = i^n \lambda_n(c) \psi_n^c (x).
\end{align}
Since $\psi^c_n(-1) \neq 0$ for all $n \in \mathbb{N}$ \cite{PSWF1}, one observes that 
\begin{align}
    f(x,y) = \frac{\psi_n^c(x)\psi_n^c(y)}{\psi_n^c(-1)} 
    \label{salut}
\end{align} verifies equation \eqref{main} and that $f(x,-1) = \psi_n^c(x)$. Thus, Theorem \ref{theo1} applies and we obtain
\begin{align}
     U(\xi ; T_x)[\psi_n^c](x) &= \left(\frac{\psi_n^c(-1 + \xi)}{\psi_n^c(-1)}\right) \psi_n^c(x).
     \label{acck}
\end{align}
This result is also a natural consequence of taking $U(\xi ; -\chi_n(c))$ in Corollary \ref{cococo}. Using the action \eqref{acck}, one finds that
\begin{align}
    \left( \int \alpha(\xi) U(\xi ; T_x) d\xi\right)[\psi_n^c] (x) =  \left(\int \alpha(\xi) \frac{\psi_n^c(-1 + \xi)}{\psi_n^c(-1)} d\xi \right) \psi_n^c(x).
    \label{prepre}
\end{align}
In particular, injecting 
\begin{align}
   \alpha(\xi) =   \left\{
    \begin{array}{ll}
    	 e^{ic(1 - \xi)}  & \mbox{if } \xi \in [0,2[,  \\
    	0 & \mbox{otherwise, }
    \end{array}
\right. 
\end{align}
in equation \eqref{prepre} yields
\begin{equation}
    \begin{split}
         \left( \int_0^2 e^{ic(1 - \xi)}  U(\xi ; T_x) d\xi\right)[\psi_n^c] (x) &=  \left(\frac{\mathcal{F}_c[\psi_n^c](-1)}{\psi_n^c(-1)} \right) \psi_n^c(x).
    \end{split}
\end{equation}
Then, it is enough to note that
\begin{align}
     \frac{\mathcal{F}_c[\psi_n^c](-1)}{\psi_n^c(-1)}  = i^n \lambda_n(c)
\end{align}
to prove the following theorem:
\begin{theorem}
Let $\mathcal{F}_c$ be the finite Fourier transform, $\mathcal{R}$ the reflection operator, $T_x$ the Heun operator defined in \eqref{heun} and $U(\xi;T_x)$ the function of $T_x$ defined by 
\begin{align}
    U(\xi ; T) \equiv \sum_{k =0}^\infty  \frac{\xi^k U_k(T_x)}{k!},
\end{align}
where the polynomials $U_k$ are given by the following four-term recurrence relation:
\begin{equation}
    \begin{split}
        U_{k+1}(T_x) &= \frac{(T_x + c^2 + k(k+1))}{2(k + 1)}U_k(T_x) \\& \quad \quad  - \frac{c^2 k}{k + 1}U_{k-1}(T_x) + \frac{ c^2 k(k-1)}{2(k + 1)} U_{k-2}(T_x),
    \end{split}
\end{equation}
and the initial condition $U_0(T_x) = 1$. Then, we have that
\begin{align}
    \mathcal{F}_c =  \int_0^2 e^{ic(1 - \xi)}  U(\xi ; T_x) d\xi
    \label{totoeqp}
\end{align}
 as an operator acting on $L^2[-1,1]$.
 \label{Toe1}
\end{theorem}
\noindent Since the series defining $U(\xi,T_x)$ converges more quickly for small $\xi$, it is interesting to note that injecting
\begin{align}
   \alpha(\xi) =   \left\{
    \begin{array}{ll}
    	 e^{ic(1 - \xi)} + \mathcal{R}e^{-ic(1 - \xi)} & \mbox{if } \xi \in [0,1],  \\
    	0 & \mbox{otherwise, }
    \end{array}
\right. 
\end{align}
in \eqref{prepre} also gives
\begin{equation}
    \begin{split}
         \left( \int_0^1  \left(e^{ic(1 - \xi)} + \mathcal{R}e^{-ic(1 - \xi)}\right) U(\xi ; T_x) d\xi\right)[\psi_n^c] (x) &=  \left(\frac{\mathcal{F}_c[\psi_n^c](-1)}{\psi_n^c(-1)} \right) \psi_n^c(x)
    \end{split}
\end{equation}
and does not use $U(\xi,T_x)$ with $\xi \in ]1,2[$. Therefore, we have
\begin{corollary}
With the same preamble as Theorem \eqref{Toe1}, we have
\begin{align}
    \mathcal{F}_c = \int_0^1  \left(e^{ic(1 - \xi)} + \mathcal{R}e^{-ic(1 - \xi)}\right) U(\xi ; T_x) d\xi
    \label{totoeq}
\end{align}
as an operator acting on $L^2[-1,1]$.
\end{corollary}
\subsection{The $\sinc$ kernel}
Given that $\mathcal{Q}_c = \frac{2\pi}{c} \mathcal{F}_c^* \circ \mathcal{F}_c$, Theorem \eqref{Toe1} is sufficient to show that $\mathcal{Q}_c$ can be expressed as a function of $T_x$. However, we would like to obtain formulas similar to \eqref{totoeqp} and \eqref{totoeq}, i.e. to find $\beta(\xi)$ such that
\begin{align}
  \mathcal{Q}_c =  \left( \int \beta(\xi) U(\xi ; T_x) d\xi\right) 
\end{align}
or equivalently 
\begin{align}
   \left( \int \beta(\xi) U(\xi ; T_x) d\xi\right)[\psi_n^c](x) = \mu_n(c)\psi_n^c(x).
\end{align}
Again, we can use \eqref{salut} and Theorem \ref{theo1} to obtain
\begin{align}
    \left( \int \beta(\xi) U(\xi ; T_x) d\xi\right)[\psi_n^c] (x) =  \left(\int \beta(\xi) \frac{\psi_n^c(-1 + \xi)}{\psi_n^c(-1)} d\xi \right) \psi_n^c(x).
    \label{prepre2}
\end{align}
Then, taking 
\begin{align}
   \beta(\xi) =   \left\{
    \begin{array}{ll}
    	 \frac{\sin{(c \xi)}}{\pi \xi} & \mbox{if } \xi \in [0,2[,  \\
    	0 & \mbox{otherwise, }
    \end{array}
\right. 
\end{align}
yields
\begin{equation}
    \begin{split}
         \left( \int_0^2 \frac{\sin{(c \xi)}}{\pi \xi}  U(\xi ; T_x) d\xi\right)[\psi_n^c] (x) &=  \left(\frac{\mathcal{Q}_c[\psi_n^c](-1)}{\psi_n^c(-1)} \right) \psi_n^c(x).
    \end{split}
\end{equation}
Since we have that
\begin{align}
    \frac{\mathcal{Q}_c[\psi_n^c](-1)}{\psi_n^c(-1)} = \mu_n(c),
\end{align}
the following theorem is proven:
\begin{theorem}
Let $\mathcal{Q}_c$ be integral operator defined in \eqref{inin}, $\mathcal{R}$ the reflection operator, $T_x$ the Heun operator defined in \eqref{heun} and $U(\xi;T_x)$ the function of $T_x$ defined by 
\begin{align}
    U(\xi ; T_x) \equiv \sum_{k =0}^\infty  \frac{\xi^k U_k(T_x)}{k!},
\end{align}
where the polynomials $U_k$ are given by the following four-term recurrence relation:
\begin{equation}
    \begin{split}
        U_{k+1}(T_x) &= \frac{(T_x + c^2 + k(k+1))}{2(k + 1)}U_k(T_x) \\& \quad \quad  - \frac{c^2 k}{k + 1}U_{k-1}(T_x) + \frac{ c^2 k(k-1)}{2(k + 1)} U_{k-2}(T_x),
    \end{split}
\end{equation}
and the initial condition $U_0(T_x) = 1$. Then, we have that
\begin{align}
    \mathcal{Q}_c = \int_0^2  \frac{\sin{(c \xi)}}{\pi \xi}  U(\xi ; T_x) d\xi
    \label{totoeq2p}
\end{align}
 as an operator acting on $L^2[-1,1]$.
 \label{Toe2}
\end{theorem}
\noindent To avoid using $U(\xi;T_x)$ for $\xi \in ]1,2[$, one could also choose
\begin{align}
   \beta(\xi) =   \left\{
    \begin{array}{ll}
    	 \frac{\sin{(c \xi)}}{\pi \xi} + \frac{\sin{(c (2- \xi))}}{\pi (2-\xi)} \mathcal{R} & \mbox{if } \xi \in [0,1],  \\
    	0 & \mbox{otherwise, }
    \end{array}
\right. 
\end{align}
to obtain
\begin{equation}
    \begin{split}
         \left( \int_0^1  \left(\frac{\sin{(c \xi)}}{\pi \xi} + \frac{\sin{(c (2- \xi))}}{\pi (2-\xi)} \mathcal{R}\right) U(\xi ; T_x) d\xi\right)[\psi_n^c] (x) &=  \left(\frac{\mathcal{Q}_c[\psi_n^c](-1)}{\psi_n^c(-1)} \right) \psi_n^c(x).
    \end{split}
\end{equation}
Then, one finds the following corollary:
\begin{corollary}
With the same preamble as Theorem \eqref{Toe2}, we have
\begin{align}
    \mathcal{Q}_c = \int_0^1  \left(\frac{\sin{(c \xi)}}{\pi \xi} + \frac{\sin{(c (2- \xi))}}{\pi (2-\xi)} \mathcal{R}\right) U(\xi ; T_x) d\xi
    \label{totoeq2}
\end{align}
as an operator acting on $L^2[-1,1]$.
\end{corollary}
\section{Limiting cases}
\label{s4}
We are now interested in cases where the formulas in Theorems \ref{Toe1} and \ref{Toe2} can be simplified. We will consider those where the bandwidth parameter $c$ is either small or large.
\subsection{ The limit $c \rightarrow 0 $}
Let us start from equation \eqref{totoeq} which can be rewritten as
\begin{align}
    \mathcal{F}_c = \int_{-1}^1  e^{-ic y} U(y + 1 ; T_x) dy.
    \label{kint}
\end{align}
From the recurrence relation \eqref{rel2}, one finds that
\begin{align}
    U_k(T_x) = \frac{\prod_{n = 1}^{k} (T_x + k(k-1))}{2^k k!} + O(c^2)
\end{align}
and thus
\begin{align}
     e^{-icy} U(y+1;T_x) = (1-icy)\sum_{k=0}^\infty \frac{\prod_{n = 1}^{k} (T_x + k(k-1))}{2^k k!k!} (y+1)^k + O(c^2).
\end{align}
Evaluating the integral in equation \eqref{kint} then yields the following:
\begin{equation}
    \begin{split}
        \mathcal{F}_c &= 2 \sum_{k=0}^\infty \frac{\prod_{n = 1}^{k} (T_x + k(k-1))}{ k!(k+1)!}\left( 1 - \frac{ick}{k+2}\right) + O(c^2)\\
        & = 2 \sum_{k=0}^\infty \frac{\prod_{n = 1}^{k} \left( (1-x^2)\pdv[2]{}{x} - 2x + k(k-1)\right)}{ k!(k+1)!}\left( 1 - \frac{ick}{k+2}\right) + O(c^2).
        \label{resc0}
    \end{split}
\end{equation}
Using Legendre polynomials $\{P_n(x)\}_{n\in \mathbb{N}}$, which give a basis of $L^2[-1,1]$ and satisfy
\begin{align}
     \left((1-x^2)\pdv[2]{}{x} - 2x \right)P_n(x) = -n(n+1)P_n(x),
\end{align}
one can check in \eqref{resc0} that $\mathcal{F}_c$ is at order $0$ in $c$ the projector onto the space of functions spanned by $P_0(x) = 1$. Similarly, the term of order $1$ in $c$ is the projector onto the space spanned by $P_1(x) = x$. Recalling the orthogonality property of the Legendre polynomials, this is indeed what is expected from the definition of $\mathcal{F}_c$ given by \eqref{fifit}:
\begin{equation}
    \begin{split}
        \mathcal{F}_c[\phi](x) &= \int_{-1}^1 \phi(y)dy  - i c \int_{-1}^1 y\phi(y)dy + O(c^2)\\
        & = \int_{-1}^1 P_0(y) \phi(y)dy  - i c \int_{-1}^1 P_1(y) \phi(y)dy + O(c^2).
    \end{split}
\end{equation}
Higher order terms in \eqref{resc0} can be obtained in a similar way.

\subsection{The limit $c \rightarrow \infty$ and the complete Fourier transform}
Let $\mathcal{D}_c$ be the dilation operator acting as:
\begin{align}
    \mathcal{D}_c \phi(x) = \phi(\sqrt{c} x)
\end{align}
When $c \rightarrow \infty$, one can check that the dilated finite Fourier transform $\Tilde{\mathcal{F}}_c = \mathcal{D}_c^{-1} \circ \mathcal{F}_c \circ \mathcal{D}_c$ yields the complete Fourier transform:
\begin{equation}
    \begin{split}
       \lim_{c \rightarrow \infty} \sqrt{ \frac{c}{2 \pi}} \Tilde{\mathcal{F}}_c[\phi](x) &= \lim_{c \rightarrow \infty} \sqrt{ \frac{c}{2 \pi}} \int_{-1}^{1} e^{i\sqrt{c} x t} \phi(\sqrt{c} t) dt \\
        &= \lim_{c \rightarrow \infty}  \frac{1}{\sqrt{2 \pi}}\int_{-c}^{c} e^{i x t} \phi(t') dt'\\
        & = \mathcal{F}[\phi](x).
    \end{split}
\end{equation}
As for the Heun operator $T_x$, under the same dilation it becomes
\begin{equation}
    \begin{split}
        \Tilde{T}_x = \mathcal{D}_c^{-1} \circ T_x \circ \mathcal{D}_c &= c\left(\pdv[2]{}{x} - x^2\right) + O(c^0) \\
        &= 2c\left(\mathcal{H} + \frac{1}{2} \right) + O(c^0).
    \end{split}
\end{equation}
Therefore, one expects that taking the limit $c \rightarrow \infty$ in equation \eqref{totoeq} should allow to recover the known result:
\begin{align}
    \mathcal{F} = e^{-i\frac{\pi}{2} \mathcal{H} },
\end{align}
where $\mathcal{H}$ is the operator defined in \eqref{hamilher}. Using equation \eqref{totoeqp}, $ y = -1 + \xi$ and conjugating by $\mathcal{D}_c$, we find:
\begin{align}
    \Tilde{\mathcal{F}}_c = \int_{-1}^{1} e^{-i c y} U(y + 1, \Tilde{T}_x) dy.
    \label{limc}
\end{align}
In the limit $c \rightarrow \infty$, we notice that the four-term recurrence relation for the polynomials $U_k$ yields
\begin{align}
    U_k =  \frac{c^{2k}}{2^k k!} + O(c^{2k -1}).
    \label{dt}
\end{align}
Then, taking $y = -1 + \epsilon/c^2$ we obtain
\begin{equation}
\begin{split}
       U(\epsilon/c^2, \Tilde{T}_x) &= \sum_{k = 0}^{\infty} \frac{\epsilon^k}{2^k k! k!}  + O(1/c)\\
& = J_0(i \sqrt{2 \epsilon}) + O(1/c),
\end{split}
\label{Bess}
\end{equation}
where $J_0$ refers to the zeroth order Bessel function.
 In particular, this expression does not depend on $\Tilde{T}_x$ and is valid as long as $y$ is near $-1$. For $\epsilon$ large, let us also note that
\begin{align}
U(\epsilon/c^2, \Tilde{T}_x) &\approx \frac{e^{\sqrt{2\epsilon}}}{\sqrt{2 \pi} (2\epsilon)^{1/4}} + O(1/c).
    \label{ee3}
\end{align}
 Outside the interval near $y = -1$, we can use Corollary \ref{cococo} to approximate $U(y+1;\Tilde{T}_x)$. As long as $y$ does not tend to $0$ or $\pm 1$, the differential equation
\begin{align}
    \Big[ (1-y^2) \pdv[2]{}{y} - 2y \pdv{}{y} - c^2 y^2 - \Tilde{T}_x \Big] U(y + 1, \Tilde{T}_x)  =0 
\end{align}
implies that
\begin{equation}
    \begin{split}
         U(y+1,\Tilde{T}_x) &= A e^{ c\sqrt{1-y^2}} \left(\frac{1}{\sqrt{y}(1-y^2)^{1/4}} \left(\frac{1+ \sqrt{1-y^2}}{1 - \sqrt{1-y^2}}\right)^{\frac{\Tilde{T}_x}{4c}} + O(1/c) \right)  \\
   &    + B e^{ -c\sqrt{1-y^2}}  \left(\frac{1}{\sqrt{y}(1-y^2)^{1/4}} \left(\frac{1+ \sqrt{1-y^2}}{1 - \sqrt{1-y^2}}\right)^{\frac{\Tilde{T}_x}{4c}}  + O(1/c) \right).
    \end{split}
        \label{ee1}
\end{equation}
The constants $A$ and $B$ are fixed by the boundary condition $U(0,T_x) = 1$. Taking $y+1 = \epsilon/c^2$ in equation \eqref{ee1}, one finds
\begin{align}
   U(y+1,T) &= A \frac{\sqrt{c} e^{ \sqrt{2\epsilon}}}{(2\epsilon)^{1/4}} \left(1 + O(1/c) \right)  \\
   &    + B \frac{\sqrt{c} e^{- \sqrt{2\epsilon}}}{(2\epsilon)^{1/4}} \left(1 + O(1/c) \right).
   \label{ee2}
\end{align}
Thus, we can compare \eqref{ee3} and \eqref{ee2} to deduce that
\begin{align}
    A = \frac{1}{\sqrt{2\pi c}}, \quad \quad B = 0.
    \label{ee4}
\end{align}
Next, we want to evaluate the integral \eqref{limc}. By introducing the complex variable $z = -iy + \sqrt{1-y^2}$, we can interpret equation \eqref{limc} as an integral
\begin{align}
    \Tilde{\mathcal{F}}_c = \int_{\Gamma_1} e^{c\frac{z^2-1}{2z}}U(1 + \frac{1-z^2}{2iz}, \Tilde{T}_x) \frac{i(1+z^2)}{2z^2} dz
\end{align}
along a path $\Gamma_1$ from $z=i$ to $z= -i$ on the half unit circle where $Re(z) > 0 $. This path can be deformed to keep away from $z =1$ ($y = 0$). This allows to use \eqref{ee1} with \eqref{ee4} to obtain
\begin{align}
    \Tilde{\mathcal{F}}_c = \int_{\Gamma_1} \frac{-e^{cz}}{\sqrt{2\pi i c}} \sqrt{\frac{1+z^2}{1-z^2}} \frac{1}{z} \left(\frac{-(1 + z)^2}{(1-z)^2}\right)^{\frac{\Tilde{T}_x}{4c}} dz
    \label{compl}
\end{align}
Along a path $\Gamma_2$ from $z = -i$ to $z = i$ in the half plane $Re(z) < 0$, the integrand tends to $0$ as $c \rightarrow \infty$ because of the term $e^{c z}$. Thus, the expression \eqref{compl} is reduced to the evaluation of its residues in the region $|z| < 1$, $Re(z) \geq 0$. Since there is only a simple pole at $z=0$, we find
\begin{align}
    \Tilde{\mathcal{F}}_c = \sqrt{\frac{2\pi i}{c}}(-1)^{\frac{\Tilde{T}_x}{4c}}(1 + O(1/c))
\end{align}
and therefore with $\Tilde{T}_x = 2c(\mathcal{H} + 1/2) + O(c^0)$, we recover
\begin{align}
    \lim_{c\rightarrow \infty} \sqrt{\frac{c}{2\pi}} \Tilde{\mathcal{F}}_c = e^{-i \frac{\pi}{2} \mathcal{H}}.
\end{align}
\begin{figure}
    \centering
\begin{tikzpicture}
\def\gap{0.2}
\def\bigradius{2}
\def\littleradius{0.3}

\draw [help lines,->] (-1.25*\bigradius, 0) -- (1.25*\bigradius,0);
\draw [help lines,->] (0, -1.25*\bigradius) -- (0, 1.25*\bigradius);
\draw [red,thick,domain=90:45,->] plot ({1.85*cos(\x)}, {1.85*sin(\x)});
\draw [red,thick,domain=46:9] plot ({1.85*cos(\x)}, {1.85*sin(\x)});
\draw [red,thick,domain=97:263] plot ({1.85 + 0.3*cos(\x)}, {0.3*sin(\x)});
\draw [red,thick,domain=270:351] plot ({1.85*cos(\x)}, {1.85*sin(\x)});
\draw [blue,thick] (-0.5, 1.85) -- (-0.25,1.85);
\draw [blue,thick] (-0.25, 1.85) -- (0,1.85);
\draw [blue,thick] (0, -1.85) -- (-0.5,-1.85);
\draw [blue,thick,->] (-0.5, -1.85) -- (-0.5,0);
\draw [blue,thick] (-0.5, 0) -- (-0.5,1.85);
\node [red] at (1.85,-1.43) {$\Gamma_{1}$};
\node [blue] at (-0.85,1.23) {$\Gamma_{2}$};

\node at (0,1.85)[circle,fill,inner sep=0.75pt]{};
\node at (0.1,2.1) {$i$};
\node at (0,-1.85)[circle,fill,inner sep=0.75pt]{};
\node at (0.2,-2.1) {$-i$};
\node at (-0.6,2.4) {$\text{Im}(z)$};
\node at (2.4,-0.3) {$\text{Re}(z)$};

\end{tikzpicture}
    \caption{Paths $\Gamma_1$ and $\Gamma_2$ in the complex plane}
    \label{fig:my_label}
\end{figure}
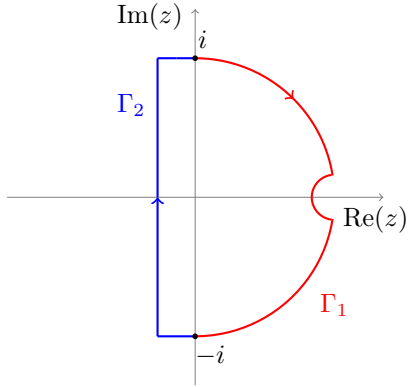

\section{Concluding remarks}

We have shown how the finite Fourier transform $\mathcal{F}_c$ can be expressed as a function of the confluent Heun operator $T$ arising in the spheroidal wave equation. In doing so, we shed new light on the relation between the two operators and have generalized the formula giving the complete Fourier transform as the exponential of a second order differential operator. 

Other settings exist in which a second order differential operator commutes with an integral one. The operator which appears in the generic Heun equation is known to commute with the finite Jacobi transform \cite{grunbaum2018algebraic}. A differential operator commuting with the finite Fourier transform for functions defined on a circle was also identified by Slepian \cite{PSWF5}. It should prove interesting to check if the approach used in this paper can be applied in those situations and if a formula relating the two commuting operators can be found. Future work could also be directed to the study of discrete cases, in which the two objects are tridiagonal matrices and complete ones \cite{grunbaum1981toeplitz}.

Finally, one expects that our results could also be derived using a more algebraic approach. Such a derivation would connect to the existing literature on the relation between the finite Fourier transform and the Heun operator, and on their associated bispectral pair in continuous and discrete settings \cite{grunbaum2018algebraic,atakishiyeva2021algebraic}.

\section*{Acknowledgments}
The authors are grateful to Alberto Grünbaum for generating their interest in time and band limiting and in the commuting operator “miracle” as he has called it. PAB holds a scholarship from the Natural Sciences and Engineering Research Council (NSERC) of Canada. NC is partially supported by Agence Nationale de la Recherche, Projet ANR-18-CE40-0001. LV gratefully acknowledges a Discovery Grant from NSERC

\end{document}